\documentclass[reqno,a4paper,12pt]{amsart}
\usepackage{amsfonts,a4}
\usepackage{amssymb} 
\usepackage{latexsym}
\usepackage{amsmath}
\usepackage{amsthm}
\usepackage{color}
\usepackage{float}

\usepackage{tikz}

\usepackage{mathrsfs}
\numberwithin{equation}{section}

\usepackage[active]{srcltx}

\newtheorem{theorem}{Theorem}[section]

\newtheorem{proposition}[theorem]{Proposition}
\newtheorem{corollary}[theorem]{Corollary}

\DeclareMathOperator{\arcsinh}{arcsinh}

\newcommand {\C}    {\mathbb{C}}

\newcommand {\R}    {\mathbb{R}}

\newcommand {\Z}    {\mathbb{Z}}

\newcommand{\DD}{{\mathbb{D}}}

\newcommand {\Ima}    {\hbox{\rm Im}\,}

\newcommand {\sign}   {\hbox{\rm sign}}

\newcommand{\el}{\hbox{\small{\rm ell}\,}}
\newcommand{\parr}{\hbox{\footnotesize{\rm par}\,}}
\newcommand{\hyp}{\hbox{\footnotesize{\rm hyp}\,}}

\renewcommand{\a}{\alpha}

\newcommand{\g}{\gamma}
\newcommand{\G}{\Gamma}

\newcommand{\e}{\epsilon}

\newcommand{\z}{\zeta}
\newcommand{\y}{\eta}

\renewcommand{\l}{\lambda}

\newcommand{\m}{\mu}
\newcommand{\n}{\nu}
\newcommand{\x}{\xi}

\renewcommand{\r}{\rho}
\newcommand{\s}{\sigma}

\newcommand{\vs}{\varsigma}
\renewcommand{\t}{\tau}

\newcommand{\h}{\chi}

\renewcommand{\u}{\upsilon}

\newcommand{\eb}{{\mathbf e}}
\newcommand{\fb}{{\mathbf f}}

\newcommand{\Eb}{{\mathbf E}}

\newcommand{\Kb}{{\mathbf K}}

\newcommand{\Mb}{{\mathbf M}}
\newcommand{\Nb}{{\mathbf N}}

\newcommand{\eF}{\mathfrak e}

\newcommand{\Ac}{{\mathcal A}}
\newcommand{\Bc}{{\mathcal B}}

\newcommand{\Dc}{{\mathcal D}}
\newcommand{\Ec}{{\mathcal E}}

\newcommand{\Jc}{{\mathcal J}}

\newcommand{\Lc}{{\mathcal L}}

\font\sc=rsfs10 at 12pt

\newcommand{\Es}{\sc\mbox{E}\hspace{1.0pt}}

\newcommand{\Ds}{\sc\mbox{D}\hspace{1.0pt}}

\newcommand{\Ss}{\sc\mbox{S}\hspace{1.0pt}}

\begin{document}

\date{\today}

\title[Commutative Toeplitz Algebras -- Spectral]{Commutative algebras of Toeplitz operators on the Bergman space revisited: Spectral theorem approach}
\author{Grigori Rozenblum}
\address{Chalmers Univ. of Technol., Sweden; The Euler Intern. Math. Institute and St.Petersburg State Univ.; Mathematics Center
Sirius Univ. of Sci. and Technol.
Sochi Russia}
\email{grigori@chalmers.se}
\thanks{G.R. was supported  by the grant of the Russian Fund of Basic Research 20-01-00451.}
\author{Nikolai Vasilevski}
\address{Department of Mathematics, CINVESTAV, Mexico City, Mexico}
\email{nvasilev@math.cinvestav.mx}
\thanks{N.V. was partially supported by CONACYT grants 280732 and FORDECYT- \\
\indent PRONACES/61517/2020238630, Mexico.}

\begin{abstract}
For three standard models of commutative algebras generated by  Toeplitz operators in the weighted analytic Bergman space on the unit disk, we find their representations as the algebras of bounded functions of certain unbounded self-adjoint operators. We discuss  main properties of these representation and, especially, describe relations between properties of  the spectral function of Toeplitz operators in the spectral representation and properties of the symbols.

\vspace{1ex}
\vspace{1ex}
\end{abstract}

\maketitle

\section{Introduction}

There exists now an advanced theory characterizing  all commutative $C^*$-algebras generated by Toeplitz operators with bounded symbols, acting on the Bergman spaces $\Ac^2_\l$ of analytic functions in the unit disk, equivalently, in the upper half-plane, see e.g. \cite{book}. All of them are classified   either by pencils of geodesics or by maximal Abelian subgroups of M\"obius transformations, which are the one-parameter groups generated by some fixed non-identical M\"obius transformation. Each of these one-parameter Abelian subgroups is conjugated to one of the following three  model groups: the \emph{elliptic} one, $\mathbb{T}$, acting on the disk $\mathbb{D}$, the \emph{parabolic} one, $\mathbb{R}$,  and the \emph{hyperbolic} one, $\mathbb{R}_+$, the latter two acting on the upper half-plane $\Pi$.
Then,  Toeplitz operators whose bounded symbols are invariant under the action of one of the above model groups, generate a commutative $C^*$-algebra, and  all commutative algebras arise as a result of this construction.

On the other hand, there exists a general abstract construction of commutative algebras of bounded operators in a Hilbert space based on the Spectral Theorem. Namely, these are the algebras of operators formed by bounded functions of some fixed self-adjoint operator.

In the paper we establish a connection between these two constructions. With  each model commutative algebra we associate a certain self-adjoint (unbounded) operator which generates this algebra  via the functional calculus. In each case this unbounded self-adjoint {\emph{generating}} operator is a first order differential operator closely related with the vector field whose integral curves are exactly the orbits of the corresponding maximal Abelian groups of M\"obious transformations.

{In the, now extensive, literature devoted to the study of commutative $C^*$-algebras generated by Toeplitz operators with symbols being invariant under the action of certain groups, the object called  the \emph{spectral function} turned out to be very useful. For each Toeplitz operator $T_a$ in  such commutative $C^*$-algebra $\mathcal{T}$,  its} spectral function is the function $\gamma_a$, such that $T_a$  is unitary equivalent to the multiplication by $\gamma_a$ acting in an appropriate $L_2$ function space: $RT_aR^*= \gamma_aI$, {where $I$ is the identity operator in this $L_2$ space.} What is important here, is that this unitary operator $R$ is the same for all $T_a$ in the algebra $\mathcal{T}$, and it depends only on the particular group defining invariant symbols, whose corresponding Toeplitz operators generate $\mathcal{T}$ and on the weight parameter $\lambda$ of the Bergman space $\Ac^2_\l$.

In \emph{this} paper we show that for each operator in these commutative $C^*$-algebras, its spectral function, thus defined, coincides with function that determines via the functional calculus this operator as a function of the (unbounded) generating operator for the corresponding algebra. This approach enables us, in particular, to extend essentially the set of symbols to rather singular ones, that nevertheless produce bounded Toeplitz operators.

Note that apart of the characterization of commutative $C^*$-algebras generated by Toeplitz operators with invariant symbols, we describe also corresponding von Neumann algebras, i.e. their strong operator topology closures, as well as  commutative algebras of unbounded invariant operators, with a natural definition of commutativity of unbounded operators.

\section{Preliminaries}
We consider the standard  weighted \emph{Bergman space} $\Ac_\l^2=\Ac_\l^2(\DD)$ of analytic functions on the unit disk $\DD\subset\C^1$ belonging to the space $L^2_\l(\DD)$ with measure $d\m_\l=\frac{\l+1}{\pi}(1-|z|^2)^{\l}dA(z),$ $\l>-1,$ where $dA$ is the Lebesgue measure. Sometimes, it is more convenient to consider the equivalent model, namely the  space $\Ac^2_\l(\Pi)$ on the upper half-plane $\Pi\subset \C^1$ with measure $d\n_\l(z)=\frac{\l+1}{\pi} (2y)^{\l}dA(z).$ This latter space is naturally isometrically isomorphic to $\Ac_\l^2(\DD).$ This isomorphism is generated by the M\"obius transformation
\begin{equation}\label{DiskToHalfplane}
  U_\l:  f(z)\mapsto \left(\frac{\sqrt{2}}{1-iz}\right)^{\l+2}f\left(\frac{z-i}{1-iz}\right).
\end{equation}

Recall, for completeness, that the standard orthogonal monomial  basis in $\mathcal{A}^2_\l(\DD)$ consists of the functions
\begin{equation}\label{basis D}
  \eb_{k}(z)\equiv  \eb_{k;\l}(z)=\sqrt{\frac{\Gamma(k+\lambda+2)}{k!\Gamma(\lambda+2)}}\,z^k.
\end{equation}

In our considerations, the weight parameter $\l>-1$ is usually  fixed; in such cases we may omit it in our notations if no confusion threatens.
By the standard definition, given a function $a \in L_{\infty}(\mathbb{D})$, the Toeplitz operator $T_a$ with \emph{symbol} $a$ is defined on $\Ac_\l^2(\DD)$ as
$T_af := B_\l(af)$, where $B_\l$ is the orthogonal Bergman projection of $L^2_{\lambda}(\mathbb{D})$ onto $\Ac_\l^2(\DD)$; the same notation is used for Toepitz operators in the half-plane model.

We recall here the known characterization of  commutative $C^*$-algebras generated by Toeplitz operators with bounded measurable symbols acting on the standard weighted Berg\-man space $\Ac_\l^2(\DD)$, see e.g. \cite{book}. As it turns out, there exists a one-to-one correspondence bet\-ween commutative $C^*$-algebras generated by Toeplitz operators and maximal Abelian subgroups of M\"obius transformations of $\DD$, which are the one-parameter groups generated by a non-identical M\"obius transformation. Each one-parameter Abelian subgroup of this kind is conjugated to one of the following three  model groups:

\noindent -  \emph{elliptic}, $\mathbb{T}: \, \mathbb{D} \rightarrow \mathbb{D}$, with action $t \in \mathbb{T}:\, z\mapsto tz$;

\noindent - \emph{parabolic}, $\mathbb{R}: \, \Pi \rightarrow \Pi$, with action $h \in \mathbb{R}:\, z\mapsto z+h$;

\noindent - \emph{hyperbolic}, $\mathbb{R}_+: \, \Pi \rightarrow \Pi$, with action $\rho \in \mathbb{R}_+:\, z\mapsto \rho z$.

The orbits of the above groups can be described, for example, by the following parametric equations:

\noindent -  \emph{group $\mathbb{T}$ on $\mathbb{D}$}: \ \ \ $x = R \cos \theta, \ y = R\sin \theta$, \ $\theta  \in [0, 2\pi)$, with a fixed $R \in (0, 1)$;

\noindent - \emph{group $\mathbb{R}$ on $\Pi$}: \ \ \ $x = h, \ y = y_0$, \ $h \in \mathbb{R}$, with a fixed $y_0 \in \mathbb{R}_+$;

\noindent - \emph{group $\mathbb{R}_+$ on $\Pi$}: \ \ { $x = \rho\cos \theta, \ y = \rho \sin \theta$, \ $\rho = e^s\in \R_+$, $s \in \mathbb{R}$, with a fixed $\theta  \in (0, 2\pi)$.}

Note that these orbits are integral curves, respectively, for the following vector fields
\begin{equation*}
 V_{\el} = -y\frac{\partial }{\partial x} + x\frac{\partial }{\partial y}, \qquad V_{\parr} = \frac{\partial }{\partial x}, \qquad V_{\hyp} = x\frac{\partial }{\partial x} + y\frac{\partial }{\partial y}.
\end{equation*}

Further \cite[Theorems 9.6.2, 10.4.1]{book}, for any maximal Abelian subgroup $G$ of M\"obius transformations, the $C^*$-algebra generated by Toeplitz operators with bounded measurable symbols which are invariant under the action of $G$ is commutative. Moreover, for each such $G$, there exist a  Borel subset $X_G \subset \mathbb{R}$ with a measure $\s_G$ and a unitary operator $R_G : \mathcal{A}^2_{\lambda} \rightarrow L_2(X_G;\s_G)$, such that each Toeplitz operator $T_a$ with $G$-invariant symbol $a$  is unitarily equivalent to the multiplication operator by a certain function $\gamma_a$, acting on $L_2(X_G;\s_G)$,
\begin{equation*}
 R_G T_a R_G^* = \gamma_aI.
\end{equation*}
This function $\g_a$ is usually (see, e.g.,  \cite{HMV})  called the \emph{spectral function} of the Toeplitz operator $T_a.$ {
We recall the description of the spectral functions for the above three model cases. In what follows we shorten the subscripts $G_{\mathrm{ell}}= \mathbb{T}$, $G_{\mathrm{par}}= \mathbb{R}$ and $G_{\mathrm{hyp}}=\mathbb{R}_+$ to $\mathrm{ell}$, $\mathrm{par}$, and $\mathrm{hyp}$, respectively.}

\noindent -  \emph{elliptic}: $X_{\el} = \mathbb{Z}_+$ with the counting measure, and (\cite[Theorem 10.3.3 and Corollary 10.3.4]{book})  the operator $R_{\el} : \, \mathcal{A}^2_{\lambda}(\mathbb{D}) \rightarrow \ell_2(\mathbb{Z}_+)$, is given by
\begin{equation*}
 R_{\el} : f(z) \ \longmapsto \ \left\{\frac{\alpha_{k,\lambda}}{\sqrt{\pi}} \int_{\mathbb{D}} f(z)\overline{z}^k d\mu_{\lambda}(z)\right\}_{k \in \mathbb{Z}_+},
\end{equation*}
where
\begin{equation*}
 \alpha_{k,\lambda} = \left( \frac{\pi \Gamma(k+\lambda+2)}{k!\Gamma(\lambda+2)} \right)^{-1/2};
\end{equation*}
its adjoint $R_{\el}^* = R_{\el}^{-1} : \, \ell_2(\mathbb{Z}_+) \rightarrow \mathcal{A}^2_{\lambda}(\mathbb{D})$ is given by
\begin{equation*}
R_{\el}^* : \  \{c_k\}_{k \in \mathbb{Z}_+} \ \longmapsto  \
\frac{1}{\sqrt{\pi}} \sum_{k \in \mathbb{Z}_+} \alpha_{k,\lambda}c_k z^k;
\end{equation*}

\noindent - \emph{parabolic}: $X_{\parr} = \mathbb{R}_+$ with the Lebesgue measure, and (\cite[Theorem 10.3.7 and Corollary 10.3.8]{book}) the operator $R_{\parr} : \, \mathcal{A}^2_{\lambda}(\Pi) \rightarrow L_2(\mathbb{R}_+)$, is given by
\begin{equation*}
 (R_{\parr}f)(\xi) = \frac{\x^{\frac{\lambda+1}{2}}}{\sqrt{\Gamma(\lambda+2)}} \int_{\Pi} f(z)e^{-i\overline{z}\xi}d\n_{\lambda}(z);
\end{equation*}
 its adjoint $R_{\parr}^* = R_{\parr}^{-1} : \, L_2(\mathbb{R}_+) \rightarrow \mathcal{A}^2_{\lambda}(\Pi)$ is given by
\begin{equation*}
(R_{\parr}^* \psi)(z) = \frac{1}{\sqrt{\Gamma(\lambda+2)}} \int_{\mathbb{R}_+} \psi(\xi) \xi^{{\frac{\lambda+1}{2}}} e^{iz \xi}d\xi;
\end{equation*}

\noindent - \emph{hyperbolic}: $X_{\hyp} = \mathbb{R}$ with the Lebesgue measure, and (\cite[Theorem 10.3.11 and Corollary 10.3.12]{book}) the operator $R_{\hyp} : \, \mathcal{A}^2_{\lambda}(\Pi) \rightarrow L_2(\mathbb{R})$, is given by
\begin{equation*}
 (R_{\hyp}f)(\eta) = \frac{\vartheta_{\lambda}(\eta)}{\sqrt{2}} \int_{\Pi} f(z)(\overline{z})^{-i\eta-(1+\lambda/2)}d\nu_{\lambda}(z);
\end{equation*}
where
\begin{equation} \label{eq:vartheta}
 \vartheta_{\lambda}(\eta) = \left(2^{\lambda}(\lambda+1) \int_0^{\pi} e^{-2\eta\theta}\sin^{\lambda}\theta d\theta\right)^{-1/2} = \frac{\left|\Gamma\left(
 \frac{\lambda+2}{2} +i\eta\right)\right|} {\sqrt{\pi\Gamma(\lambda+2)}}e^{\pi\eta/2},
\end{equation}
its adjoint $R_{\hyp}^* = R_{\hyp}^{-1} : \, L_2(\mathbb{R}) \rightarrow \mathcal{A}^2_{\lambda}(\Pi)$ is given by
\begin{equation*}
(R_{\hyp}^* g)(z) = \frac{1}{\sqrt{2}} \int_{\mathbb{R}} g(\eta) z^{i\eta-(1+\lambda/2)} \vartheta_{\lambda}(\eta) d\eta.
\end{equation*}

The corresponding spectral functions $\gamma_a$ for invariant symbols $a$ are given (\cite[Chapter 10]{book}) respectively by

\noindent -  \emph{elliptic}: $a=a(|z|)$, \ $z  \in \mathbb{D}$, \ \ $\gamma_a = \{\gamma_a(k)\}_{k \in \mathbb{Z}_+}$, where
\begin{equation}\label{gamma ell}
 \gamma_a(k) = \frac{\Gamma(k+\lambda+2)}{\Gamma(k+1) \Gamma(\lambda+1)} \int_0^1 a(\sqrt{r})r^k(1-r)^{\lambda}dr;
 \end{equation}

\noindent -  \emph{parabolic}: $a=a(\mathrm{Im}\,z)$, \ $z  \in \Pi$,
\begin{equation}\label{gamma par}
 \gamma_a(\eta) = \frac{1}{\sqrt{\Gamma(\lambda+1)}} \int_{\mathbb{R}_+} a\left(\frac{y}{2\eta}\right) y^{\lambda} e^{-y}dy,\qquad \eta \in \mathbb{R}_+;
 \end{equation}

\noindent -  \emph{hyperbolic}: $a=a(\mathrm{arg}\,z) = a(\theta)$,  \ $z = |z|e^{i\theta} \in \Pi$,
\begin{equation}\label{gamma hyp}
 \gamma_a(\eta) = 2^{\lambda}(\lambda+1)\vartheta_{\lambda}^2(\eta) \int_0^{\pi} a(\theta) e^{-2\eta\theta} \sin^{\lambda}\theta \,d\theta, \qquad \eta \in \mathbb{R};
\end{equation}
Here, due to the known asymptotics of the $\G$-function (see, e.g., \cite{Erd}, Sect. 1.18 (2),)
\begin{equation}\label{Gamma}
 \vartheta_{\lambda}(\eta)^2=O(e^{\pi\eta-(\pi-\e)|\eta|}),\, |\eta|\to \infty,
\end{equation}
for an arbitrary $\e>0.$
The latter relation is used to establish estimates for $\gamma_a(\eta)$ in concrete situations.
Defined initially { for Toeplitz operators with bounded symbols $a$ possessing the corresponding invariance properties, the mapping
\begin{equation*}
    \Ss\,: T_a\mapsto \g_a
\end{equation*}
is naturally extended to Toeplitz operators with wider} classes of symbols $a,$ namely, with locally bounded ones but unbounded near the boundary, and, further on, with the distributional ones, so that the spectral function $\g_a$ is still bounded, provided proper modifications of formulas \eqref{gamma ell}, \eqref{gamma par}, \eqref{gamma hyp} are made. This procedure expands the set of symbols for which the corresponding Toeplitz operators are bounded.

The image   of the set of bounded functions $\g_a$ under the mapping $\Ss$ is hard to describe due to many reasons. In particular, the mapping $\Ss$ does not respect the multiplication, { $\Ss(T_{a_1}T_{a_2}),$ is generally  not equal to $\Ss(T_{a_1})\Ss(T_{a_2}).$}
One can, however, { extend the mapping $\Ss$ to all operators in the $C^*$-algebra generated by Toeplitz operators with bounded invariant symbols.
In this setting the explicit characterization of the image of the mapping $\Ss$ is known for all three model cases of the above commutative $C^*$-algebras (see the forthcoming sections).}

\section{The elliptic case, radial operators}

This case has been already  considered in \cite[Subsection 4.6]{Vas_prerp-22}, where basic proofs and details can be found. We present here a further discussion with more insight.
We start with the unbounded Hermitian operator $\frac{1}{i}V_{\el} = z \frac{\partial }{\partial z} - \overline{z}\frac{\partial }{\partial \overline{z}} = \frac{1}{i} \frac{\partial }{\partial \theta}$ in $L_2(\mathbb{D}, d\mu_{\lambda})$, introduced in \cite[Formula $(4.9)$]{Wunsche},  and called there the operator of the (orbital) angular momentum. This operator, defined initially on smooth functions in $L_2(\mathbb{D}, d\mu_{\lambda}),$ is known to be closable and essentially self-adjoint. We keep the notation $\frac{1}{i}V_{\el}$ for the closure of $\frac{1}{i}V_{\el}$.

 The standard polynomial basis $e_{k,l}(z),$ $k,l\ge 0,$ in $L_2(\mathbb{D}, d\mu_{\lambda}),$ is described in detail in \cite{Wunsche}, see also \cite[Formula $(32)$]{BMR}; these functions are expressed via the shifted Jacobi polynomials and are obtained by means of the orthonormalization of the monomial system $z^k \bar{z}^l$. The explicit expression for these functions is of no importance here, with the exception for the ones with $l=0.$ Namely, we have $e_{k,0}(z)=\eb_k(z)$ {given by \eqref{basis D} and forming the standard monomial basis in the Bergman space $\Ac^2_\l(\DD).$}

The operator $\frac{1}{i}V_{\el}$ is well-defined on the linear span of the basis elements of $L_2(\mathbb{D}, d\mu_{\lambda})$ and acts on them as 
\begin{equation*}
 \textstyle{\frac{1}{i}}V_{\el} e_{k,l}: = (k-l) e_{k,l},\qquad \mathrm{for \ all} \quad (k,l) \in \mathbb{Z}_+\times \mathbb{Z}_+.
\end{equation*}
Then the operator $\frac{1}{i}V_{\el}$ can be extended to the self-adjoint operator on the domain
\begin{equation*}
 \Dc(\textstyle{\frac{1}{i}}V_{\el}) = \left\{\fb = \sum_{k,l \in \mathbb{Z}_+}\fb_{k,l}e_{k,l} \in L_2(\mathbb{D}, d\mu_{\lambda}) : \, \sum_{k,l \in \mathbb{Z}_+\times  \mathbb{Z}_+}|(k-l) \fb_{k,l}|^2 < \infty \right\}
\end{equation*}
by the same rule
\begin{equation}\label{momentum}
 \textstyle{\frac{1}{i}}V_{\el}\fb = \sum_{k,l \in \mathbb{Z}_+\times  \mathbb{Z}_+}(k-l)\fb_{k,l}e_{k,l}.
\end{equation}
Defined in this way, the self-adjoint operator $\frac{1}{i}V_{\el}$ coincides with the one constructed above by closing  the operator defined on smooth functions.

The operator $\frac{1}{i}V_{\el}$  commutes obviously with the rotation operators $\varpi(t)\fb(z) = \fb(t^{-1}z)$, $t \in \mathbb{T}$, providing thus an example of an unbounded \emph{radial} operator.

It follows from   representation \eqref{momentum} that
\begin{equation*}
    \textstyle{\frac{1}{i}}V_{\el}: \ \Dc\left(\textstyle{\frac{1}{i}}V_{\el}\right)\cap\mathcal{A}^2_{\lambda}(\mathbb{D})\to \mathcal{A}^2_{\lambda}(\mathbb{D}),
\end{equation*}
therefore, we can consider the restriction $\Nb$ of $\frac{1}{i}V_{\el}$ to $\Dc(\Nb):=\Dc\left(\frac{1}{i}V_{\el}\right)\cap\mathcal{A}^2_{\lambda}(\mathbb{D}).$
 This restriction is a self-adjoint operator in $\mathcal{A}^2_{\lambda}(\mathbb{D}),$
whose action  on  $\Dc(\Nb)$ is described by
\begin{equation}\label{N}
 \Nb\left(\sum_{k \in \mathbb{Z}_+} \fb_{k}e_{k,0}\right) \equiv \Nb\left(\sum_{k \in \mathbb{Z}_+} \fb_{k}\eb_{k}\right) = \sum_{k \in \mathbb{Z}_+} k\,\fb_{k}\eb_{k}.
\end{equation}
Since $\frac{\partial }{\partial \overline{z}}=0$ on $\mathcal{A}^2_{\lambda}(\mathbb{D}),$ the action of the operator $\Nb$ can be also described as $(\Nb f)(z)= z\frac{\partial }{\partial {z}}f(z)$, and in polar co-ordinates $r,\theta$ the operator $\Nb$ acts as $\Nb : \,f(r,\theta)\mapsto \frac{1}{i}\frac{\partial }{\partial {\theta}}f(r,\theta)$.

As follows from \eqref{N}, the system $\eb_{k},$ $k\ge0,$ diagonalizes the operator $\Nb.$ Thus, the spectrum of $\Nb$ consists of the eigenvalues $k\in \Z_+$ with corresponding eigenfunctions $\eb_{k};$ all eigenvalues are simple. We denote by $\mathrm{P}_k$ the spectral projection of $\Nb$ corresponding to the eigenvalue $k;$ it is a rank-one projection, $\mathrm{P}_k f =\langle f, \eb_k\rangle \eb_k.$

The \emph{spectral measure} $E=E(S):=E_{\Nb}(S)$, where $S$ is a Borel subset in $ \mathbb{R}^1$, of the operator $\Nb$ is given by
\begin{equation*}
 E(S) = \sum_{k \in S} \mathrm{P}_k.
\end{equation*}
So, the set $(-\infty,0) \bigcup_{k \in \Z_+} (k, k+1)$ has zero $E$-measure, while $E(\{k\}) = \mathrm{P}_k$.   For the \emph{resolution of identity} of $\Nb$, $\Eb(\y):= E((-\infty,\y)),$ we have $\Eb(\y)=\sum_{k<\y}\mathrm{P}_k.$

By the spectral theorem, we can write
\begin{equation*}
 \Nb = \int_{\mathbb{R}}\eta\, d\Eb(\eta) = \sum_{k \in \Z_+} k\, \mathrm{P}_k,
\end{equation*}
the integral understood in the sense of the strong convergence. Further, for each $f = \sum_{k \in \Z_+}f_k\eb_{k} \in \mathcal{A}^2_{\lambda}(\DD)$, we introduce the bounded, non-decreasing left continuous step function
\begin{equation*} 
 \rho_f(\eta) := \langle f, \Eb(\eta)f\rangle = \| \Eb(\eta)f\|^2 = \sum_{k < \eta}|f_k|^2.
\end{equation*}
This function $\rho_f(\eta)$ defines in the standard way  the Stieltjes measure on $\mathbb{R}$, which we will also denote by $\rho_f$. Recall  that a (complex-valued Borel) function $\varphi: \mathbb{R} \rightarrow \mathbb{C}$ is called $E$-measurable if it is $\rho_f$ - measurable for all $f \in \mathcal{A}^2_{\lambda}(\mathbb{D})$.\\
Note that for each $f = \sum_{k \in \Z_+}f_k\eb_k \in \mathcal{A}^2_{\lambda}(\mathbb{D})$, the set $(-\infty,0) \bigcup_{k \in \Z_+} (k, k+1)$ has zero $\rho_f$-measure, while $\rho_f(\{k\}) = |f_k|^2$.

It follows that, after identifying $E$-measurable functions $\varphi$ that differ on a zero $E$-measure set,  each such class of equivalent functions is uniquely defined by the sequence $\pmb{\varphi}= \{\varphi(k)\}_ {k \in \Z_+}$.

The Spectral Theorem implies now

\begin{proposition}\label{Prop.Funct.Calc.Rad}
 Given a $E$-measurable function $\varphi$, the operator
 \begin{equation*}
  \varphi(\Nb) = \int_{\mathbb{R}} \varphi(\eta)\,d\Eb(\eta) = \sum_{k \in \Z_+} \varphi(k)\mathrm{P}_k
 \end{equation*}
is well defined and normal on its domain
\begin{equation*}
 \Dc_{\varphi} = \left\{f = \sum_{k \in \Z_+}f_k\eb_k \in \mathcal{A}^2_{\lambda}(\DD) : \int_{\mathbb{R}} |\varphi(\eta)|^2\, d\rho_f(\eta) = \sum_{k \in \Z_+} |\varphi(k)f_k|^2 < \infty\right\}.
\end{equation*}
The operator $\varphi(\Nb)$ is bounded, and thus defined on the whole $\mathcal{A}^2_{\lambda}(\DD)$, if and only if the sequence   $\varphi= \{\varphi(k)\}_ {k \in \Z_+}$ is bounded.
\end{proposition}

We can now apply Theorem 1 in Sect.5.3 in \cite{BirSol}.

\begin{proposition}The mapping $\Jc: \varphi\mapsto \Jc(\varphi):=\varphi(\Nb)$ is an isometric isomorphism of the unital $C^*$-algebra $L_\infty(\R,E)$ with involution $\varphi\mapsto \bar{\varphi}$ onto a unital commutative subalgebra in $\Bc(\Ac_\l^2(\DD))$ with involution $H\mapsto H^*$.
\end{proposition}

\begin{corollary}
 The set of all operators $\{\varphi(\Nb)\}$, defined by (the classes of equivalency) of  $E$-measurable   functions  $\varphi$ with $\varphi= \{\varphi(k)\}_ {k \in \mathcal{Z}_+} \in \ell_{\infty}$ constitutes the class of  bounded mutually commuting radial operators in $\mathcal{A}^2_{\lambda}(\mathbb{D})$,  coinciding thus  with the von Neumann algebra of bounded radial operators in $\mathcal{A}^2_{\lambda}(\mathbb{D})$.
\end{corollary}

\smallskip
Let us make several remarks on the structure of the von Neumann algebra $\Lc=\Jc(\ell_{\infty}).$
First of all, recall that any \emph{bounded} radial   operator in $\Ac_\l^2(\DD)$ is diagonal in the basis $\eb_k, \, k\in \Z_+,$ therefore it has the form $\Jc(\varphi)$ for a certain $\varphi\in \ell_{\infty}.$

In particular, the Toeplitz operator with bounded radial symbol $a(|z|)\in L_{\infty}(\DD)$  (thus $a =  a(r) \in L_{\infty}(0,1)$) is the image under the mapping $\Jc$ of the function $\varphi_a\in L_\infty(\mathbb{R},E)$, given by
\begin{equation*} 
  \varphi_a(\eta) =
  \begin{cases}
   0,  & \mathrm{for} \quad \eta <0 \\ {\displaystyle
   \frac{\Gamma(\eta+\lambda+2)}{\Gamma(\eta+1) \Gamma(\lambda+1)}\int_0^1} a(\sqrt{r})r^{\eta}(1-r)^{\lambda}dr,  & \mathrm{for} \quad \eta \geq 0
  \end{cases},
 \end{equation*}
and therefore belongs to $\Lc.$

 The sequence $\varphi_a(k), k\in \Z_+,$ is, obviously, bounded and the radial Toeplitz operator $T_a$, with symbol $a$, admits the representation
\begin{equation*}
 T_a = \varphi_a(\Nb) = \int_{\mathbb{R}} \varphi_a(\eta)\,dE(\eta) = \sum_{k \in \mathcal{Z}_+} \varphi_a(k)\mathrm{P}_k,
\end{equation*}
where $\varphi_a(k)$ is given by
\begin{equation}\label{varphi_ell}
 \varphi_a(k) = \frac{\Gamma(k+\lambda+2)}{\Gamma(k+1) \Gamma(\lambda+1)} \int_0^1 a(\sqrt{r})r^k(1-r)^{\lambda}dr,
 \end{equation}
 and coincides with $\gamma_a(k)$ in \eqref{gamma ell}.

 On the other hand, such Toeplitz operators do not exhaust $\Lc.$ Namely, for a radial Toeplitz operator to be bounded, boundedness of the symbol $a(r)$ is not necessary, and this symbol may be rather singular. In fact, a, properly defined, Toeplitz  operator with radial symbol being a distribution with compact support inside $\DD$ is bounded (see  more general considerations in \cite{RozVas}.)
 Therefore, if $a(z)=a(|z|)\in \Es'(\DD)$ is a radial distribution, the expression \eqref{varphi_ell} for  the  corresponding function  $\varphi_a(k)$ providing $T_a =\varphi_a(\Nb)$ can be rewritten as
 \begin{equation*}
    \varphi_a(k)=\frac{\Gamma(k+\lambda+2)}{\Gamma(k+1) \Gamma(\lambda+1)}\pmb{(}a,|z|^{2k}(1-|z|)^\l\pmb{)},
 \end{equation*}
 where the parentheses $\pmb{(}.,.\pmb{)}$ denote the action of a distribution in $\Es'(\DD)$ on a smooth function, and this function is bounded. Even for a distribution $a$ in larger class $\Ds'(\DD),$ with the compact support condition dropped, i.e., with support touching the boundary of the disk $\DD,$ the Toeplitz operator with symbol $a,$ properly defined by means of the sesquilinear form, can be bounded, see \cite{RozVas} (where the case of $\l=0$ is considered, but the results are easily carried over to general weights.)

 Further, even among locally bounded but \emph{unbounded} (as the point $z$ approaches $\partial \DD$)
 radial symbols, there also exist examples of those that generate bounded operators.  By \cite[Example 6.1.7]{book}, the Toeplitz operators $T_a$ with unbounded but rapidly oscillating symbol
\begin{equation*}
 a(r) = (1 - r^2 )^{-\beta} \sin(1 - r^2 )^{-\alpha}, \qquad \mathrm{with} \quad 0 < \beta < \alpha,
\end{equation*}
is compact, since its eigenvalue sequence $\varphi_a(k)$  belongs to $c_0$.
Generally, bounded radial Toeplitz operators with this kind of the boundary behavior of  symbols may but need not belong to
the $C^*$-algebra generated by Toeplitz operator with \emph{bounded} radial symbols. Considerations in \cite[Proposition 5.2]{GMV} provide an example of \emph{unbounded} radial symbols generating \emph{bounded} radial Toeplitz operators which do not belong to the $C^*$-algebra generated by Toeplitz operator with bounded radial symbols.

Returning to Toeplitz operators with bounded radial symbols, it is possible to give a complete description of the those $\varphi$ for which
the operator $\varphi(\Nb)$ belongs to the $C^*$-algebra generated by Toeplitz operator with such symbols. This may happen (see e.g. \cite{GMV}) if and only if the sequence $\varphi= \{\varphi(k)\}_ {k \in \mathcal{Z}_+}$ belongs to  $ \mathrm{SO}(\mathbb{Z}_+)$, the set of all bounded sequences
which \emph{slowly oscillate} in the sense of
Schmidt \cite{Schmidt_1924_divergente_Folgen}
\begin{equation*}
 \mathrm{SO}(\mathbb{Z}_+) =
\Bigl\{\pmb{\vs}\in\ell_\infty\colon\ \lim_{\frac{j+1}{k+1}\to1}|\vs_j-\vs_k|=0\Bigr\}.
\end{equation*}
Since $\mathrm{SO}(\mathbb{Z}_+)$ is a proper subset of $\ell_{\infty}$, there exist bounded radial operators that do not belong to the $C^*$-algebra generated by Toeplitz operator with bounded radial symbols. A simple example is provided by the reflection operator $(\mathbf{R} f)(z) = f(-z)$. For this operator,  $\mathbf{R}=\varphi(\Nb)$ with function $\varphi(k)= (-1)^{k};$ such sequence belongs to $\ell_{\infty}$, but does not belong to $\mathrm{SO}(\mathbb{Z}_+)$. Moreover, this operator cannot even be expressed as a Toeplitz operator with  distributional symbol in $\Ds'(\DD)$.

 If we drop the condition of the essential boundedness of the function $\varphi$, so that the sequence $\{\varphi(k)\}_{k \in \mathbb{Z}_+}$ may be unbounded, we obtain a mapping of the space of equivalence classes of $E$-measurable functions to the set of \emph{unbounded} normal radial operators in $\mathcal{A}^2_{\lambda}(\mathbb{D})$. These operators  commute being restricted to the invariant  lineal of finite linear combinations of functions $\eb_k.$ Moreover, those operators $\varphi(\Nb)$ for which there exists at least one regular point (this means that the set of complex numbers $\varphi(k)\, , k\in \Z_+$ is not dense in $\C$) are commuting in the resolvent sense{, i.e. bounded on $\mathcal{A}^2_{\lambda}(\mathbb{D})$ resolvents in corresponding regular points commute.} {Note that the condition for the existence of at least one regular point cannot be excluded}: it is quite possible that the set $\{\varphi(k)\}_{k \in \mathbb{Z}_+}$ is dense in $\C^1$. Just numerate,
  $z_1,z_2,\dots,$ all complex numbers with both real and imaginary parts rational and set $\varphi(k)=z_k.$

Finally, as known, the radial Toeplitz operator $T_a$ is compact if and only if the sequence $\varphi(k)$ tends to zero as $k\to\infty.$

\section{The parabolic case, vertical operators}

We start with the unbounded Hermitian operator $\frac{1}{i} \frac{\partial }{\partial x} = \frac{1}{i}V_{\parr} = \frac{1}{i} \left(\frac{\partial }{\partial z} + \frac{\partial }{\partial \overline{z}}\right)$ densely defined in $L_2(\Pi, d\nu_{\lambda})\cap C^1(\Pi)$. It commutes with all translation operators $\tau_h f(z) = f(z-h),\, h\in\R$, \ $f\in  L_2(\Pi, d\nu_{\lambda}),$ providing thus an example of an unbounded vertical operator.

We restrict the operator $\frac{1}{i}V_{\parr}$ to the Bergman space $\mathcal{A}^2_{\lambda}(\Pi)$ and denote this restriction by $\Kb;$ this operator acts as  $\Kb = \frac{1}{i}\frac{\partial }{\partial z},$ considered on the natural domain  $\Dc(\Kb)=\{f\in\mathcal{A}^2_{\lambda}(\Pi): \, \Kb f\in \mathcal{A}^2_{\lambda}(\Pi) \}.$ We can observe that
\begin{gather*}
 \left(R_{\parr}\frac{1}{i}\frac{\partial }{\partial z}R_{\parr}^*\right) \psi = R_{\parr}\left(\frac{1}{i}\frac{\partial }{\partial z}\frac{1}{\sqrt{\Gamma(\lambda+2)}} \int_{\mathbb{R}_+} \psi(\xi) \xi^{{\frac{\lambda+1}{2}}} e^{iz \xi}d\xi\right) \\\nonumber = (R_{\parr}  R_{\parr}^*)(\xi\psi) = \xi \psi,\, \psi\in L_2(\R_+).
 \end{gather*}
This relation implies that the operator $\Kb = \frac{1}{i}\frac{\partial }{\partial z}$ defined on the domain $\Dc(\Kb)$
is, in fact,  unitary equivalent to the self-adjoint operator $\xi I$, $\frac{1}{i}\frac{\partial }{\partial z} = R_{\parr}^* \xi R_{\parr}$ and, therefore  self-adjoint as well.

It follows that $\Kb$ has purely absolutely continuous spectrum { of multiplicity 1}, filling the positive semi-axis $\R_+$.  The spectral measure $E=E_\Kb$ of $\Kb$ is given by
\begin{equation}\label{Sp.Measure K}
    (E(S)f)(z)=((R_{\parr})^*\h_S(\x) R_{\parr} f)(z),
\end{equation}
where $\h_S(\x)$ is the characteristic function of the Borel set $S\subset\R_+.$
Using the explicit expression for the unitaries $R_{\parr},$ $R_{\parr}^*$,
 we obtain the representation for $E(S):$
 \begin{equation*} 
  (E(S)f)(z)=\frac{1}{\sqrt{\Gamma(\lambda+2)}} \int_{S} \xi^{{\frac{\lambda+1}{2}}} e^{iz \xi}d\xi\frac{\x^{\frac{\lambda+1}{2}}}{\sqrt{\Gamma(\lambda+2)}} \int_{\Pi} f(\z)e^{-i\overline{\z}\xi}d\n_{\lambda}(\z).
 \end{equation*}
Therefore, $E(S)$ is an integral operator,
\begin{equation}\label{Sp.Function K.2}
 (E(S)f)(z)=\int_{\Pi}\eF_S(z,\z)f(\z)d\n_{\lambda}(\z),
\end{equation}
where the (distributional) kernel $\eF_S(z,\z)$ is given by
\begin{equation*}
 \eF_S(z,\z)= \frac{1}{\sqrt{\Gamma(\lambda+2)}} \int_{S}  \xi^{{\frac{\lambda+1}{2}}} e^{iz \xi} \frac{\x^{\frac{\lambda+1}{2}}}{\sqrt{\Gamma(\lambda+2)}} e^{-i\bar{\z}}\xi^{{\frac{\lambda+1}{2}}}d\xi =
\frac{1}{\Gamma(\lambda+2)}\int_{S}\xi^{\lambda+1}e^{i(z-\overline{\z})\x}d\x.
\end{equation*}

Being the spectral measure of the self-adjoint operator $\Kb,$ $E(S)$ is an orthogonal projection, which, of course, can be seen directly from \eqref{Sp.Function K.2} with some effort. A more explicit expression can be obtained for the resolution of the identity for $\Kb:$
\begin{equation*} 
 (\Eb(\eta)f)(z):=(E((0,\eta))f)(z)= \int_{\Pi}\Ec_{\eta}(z,\z)f(\z)d\n_{\lambda}(\z),
\end{equation*}
with kernel
\begin{equation*}
 \Ec_{\y}(z,\z)=\frac{1}{\Gamma(\lambda+2)}\int_0^\eta \xi^{\lambda+1}e^{i(z-\overline{\z})\x}d\x.
\end{equation*}
This kernel can be expressed via the incomplete Gamma function
\begin{equation*}
    \g(\a,s)=\int_0^s\t^{\a-1}e^{-\t}d\t
\end{equation*}
in the following way. Denote $z-\bar{\z}$ by $i\u,$ $\mathrm{Re}\,\u>0.$ Then
\begin{eqnarray} \label{EF.1}
 \Ec_{\y}(z,\z) &=&   \frac{1}{\Gamma(\lambda+2)}\int_0^\y \xi^{\lambda+1}e^{-\u\x}d\x=\
 \frac{1}{\Gamma(\lambda+2)}\u^{-(\l+2)}\int_0^{\u^{-1}\y}s^{\l+1}e^{-s}ds \\  \nonumber
 &=& \frac{1}{\Gamma(\lambda+2)}\u^{-(\l+2)}\g(\l+2,\u^{-1}\y).
\end{eqnarray}

We recall the modified incomplete Gamma function, $\g^*(\a,\varrho)=\varrho^{-\a}\G(\a)^{-1}\g(\a,\varrho)$, which is an entire analytic
function of both variables (see \cite{Erd}); taken this into account, \eqref{EF.1} gives
\begin{equation}\label{EF.2}
\Ec_{\y}(z,\z)= \y^{\l+2}\g^*\left(\l+2, \frac{\y}{i(z-\bar{\z})} \right).
\end{equation}

The Spectral Theorem leads to

\begin{proposition}
 Given a $E$-measurable function $\varphi$, the operator
 \begin{equation*}
  \varphi(\Kb) = \int_{\mathbb{R}} \varphi(\eta)\,d\Eb(\eta)
 \end{equation*}
is well defined and normal on its domain
\begin{equation*}
 \Dc_{\varphi} = \left\{f \in \mathcal{A}^2_{\lambda}(\Pi) : \int_{\mathbb{R}} |\varphi(\eta)|^2\, d\rho_f(\eta)  < \infty\right\},
 \end{equation*}
 where the measure $d\rho_f(\eta)$ is defined by the function
 \begin{equation*} 
 \rho_f(\eta) := \langle f, \Eb(\eta)f\rangle = \| \Eb(\eta)f\|^2.
\end{equation*}

The operator $\varphi(\Kb)$ is bounded, and thus defined on the whole $\mathcal{A}^2_{\lambda}(\Pi)$, if and only if the function   $\varphi$ is $E$-essentially  bounded.
\end{proposition}

Again, Theorem 1 in Sect.5.3 in \cite{BirSol} implies

\begin{proposition}The mapping $\Jc: \varphi\mapsto \Jc(\varphi):=\varphi(\Kb)$ is an isometric isomorphism of the unital $C^*$-algebra $L_\infty(\R,E)$ with involution $\varphi\mapsto \bar{\varphi}$ onto a commutative unital subalgebra in $\Bc(\Ac_\l^2(\Pi))$ with involution $H\mapsto H^*$.
\end{proposition}

\begin{corollary}
 The set of all operators $\{\varphi(\Kb)\}$, defined by (the classes of equivalency) of  $E$-measurable   functions  $\varphi \in L_\infty(\R,E)$ constitutes the class of  bounded mutually commuting vertical operators in $\mathcal{A}^2_{\lambda}(\Pi)$, coinciding thus
 with the von Neumann algebra $\Jc(L_\infty(\R,E))$ of bounded vertical operators in $\mathcal{A}^2_{\lambda}(\Pi)$.
\end{corollary}

Recall that any \emph{bounded} vertical Toeplitz operator in $\Ac^\l_2(\Pi)$ is unitary equivalent to the multiplication operator in $L_2(\R_+)$, having therefore the form $\Jc(\varphi)$ for a certain $\varphi\in L_\infty(\R,E).$

In particular, the Toeplitz operator with bounded measurable vertical symbol $a(\Ima z)$ is the image under the mapping $\Jc$ of the function $\varphi_a\in L_\infty(\R,E)$, given by
\begin{equation*}
  \varphi_a(\eta)=
  \begin{cases}
   0,  & \mathrm{for} \quad \eta <0 \\
   \frac{1}{\sqrt{\Gamma(\lambda+1)}} {\displaystyle\int_{\mathbb{R}_+}} a\left(\frac{y}{2\eta}\right) y^{\lambda} e^{-y}dy, & \mathrm{for} \quad \eta \ge 0
  \end{cases},
 \end{equation*}
and therefore belongs to $\Lc=\Jc(L_\infty(\R,E))$.

Note that, in case the function $\varphi_a$ is bounded, the vertical Toeplitz operator $T_a$, with symbol $a=a(\Ima z)$, admits the representation
\begin{equation*}
 T_a = \varphi_a(\Kb) = \int_{\mathbb{R}} \varphi_a(\eta)\,d\Eb(\eta) = R_{\parr}^*\varphi_a R_{\parr},
\end{equation*}
where $\varphi_a(\eta)$ is given by
\begin{gather}\label{eq:phi_a}
 \varphi_a(\eta) = \frac{1}{\sqrt{\Gamma(\lambda+1)}} \int_{\mathbb{R}_+} a\left(\frac{y}{2\eta}\right) y^{\lambda} e^{-y}dy\\\nonumber
 =\frac{(2\y)^{\l+1}}{\sqrt{\Gamma(\lambda+1)}}\int_{\mathbb{R}_+}a(s)s^{\l}e^{-2\y s}ds, \, \y\ge0.
 \end{gather}
This means that $\varphi_a(\eta)$ coincides with $\gamma_a(\eta)$ in \eqref{gamma par}.

At the same time, such Toeplitz operators do not exhaust $\Lc$. There exist unbounded symbols $a=a(\Ima z)$ for which the corresponding vertical Toeplitz operators are bounded, see e.g. \cite[Section 6]{HMV} for the corresponding example.

 Moreover, a wide class of bounded Toeplitz operators is obtained by considering \emph{distributional} vertical symbols, see \cite{RVhp},\cite{ERV} where such operators have been discussed both for $a(y)\in \Es(\R_+)$ and for a wide class of $a(y)\in \Ds'(\R_+)$. For example, for a compactly supported distribution $a\in \Es'(\R_+),$ consider the distribution $a^{\dag}=\pmb{1}\otimes a.$  By \eqref{eq:phi_a}, the spectral function $\varphi_{a^\dag}(\y)$ is bounded. Namely, as known, see, e.g., Theorem in Sect.4.4 Ch.II., \cite{G2}, since $a$ has compact support in $\R_+$, there exist a finite collection of functions $b_l\in C_0(\R_+)$ and constant coefficients differential operators $D_l$  such that $a=\sum_l D_l b_l, $ with derivative understood in the sense of distributions. Therefore,
 \begin{equation} \label{a distrib}
  \varphi_a(\eta)= \frac{(2\y)^{\l+1}}{\sqrt{\Gamma(\lambda+1)}}\pmb{(}\sum_l D_l b_l, s^{\l}e^{-2\eta s}\pmb{)}=
   \frac{(2\y)^{\l+1}}{\sqrt{\Gamma(\lambda+1)}}\sum_l\pmb{(} b_l, D_l^*s^{\l}e^{-2\eta s} \pmb{)},
 \end{equation}
and this is a bounded function. This implies, in particular, that the operator $\varphi_a(\Kb)$ is bounded.

As an example, we consider $a = \delta(y-1) = \frac{1}{2}\frac{d }{d y} \sign(y-1)$. Then, by \eqref{a distrib},
\begin{equation} \label{eq:varphi_delta}
 \varphi_a(\eta)= \frac{(2\y)^{\l+1}}{\sqrt{\Gamma(\lambda+1)}}\pmb{(}\delta(s-1), s^{\l}e^{-2\eta s}\pmb{)} = \frac{(2\y)^{\l+1}}{\sqrt{\Gamma(\lambda+1)}} e^{-2\eta}.
\end{equation}

{Note that the compact support condition for the distribution $a$ can be considerably relaxed, see \cite{RVhp}, \cite{ERV}.} In particular, in the above examples the distributional symbols, generating bounded Toeplitz operators, have compact support only in the vertical direction in the upper half-plane,  becoming non-compactly supported in the disk after the M\"obius transformation \eqref{DiskToHalfplane}. This produces a new set of examples of bounded (but not compact!) Toeplitz operators  with distributional symbols on the disk, not covered by the considerations in \cite{RozVas}.

\smallskip
Returning to Toeplitz operators with bounded vertical symbols, the complete description of the those $\varphi$ for which
the operator $\varphi(\Kb)$ belongs to the $C^*$-algebra generated by Toeplitz operator with such symbols is given in \cite{HMV}. This happens (see e.g. \cite{GMV}) if and only if $\varphi$ belongs to  $\mathrm{VSO}(\mathbb{R}_+)$, the set of all bounded functions
which satisfy the condition
\begin{equation*}
 \mathrm{VSO}(\mathbb{R}_+) =
\Bigl\{\varphi\in L_\infty(\R_+) \colon\ \lim_{\frac{t'}{t}\to1}|\varphi(t)-\varphi(t')|=0\Bigr\},
\end{equation*}
or, equivalently, those functions $\varphi$ that are uniformly continuous on $\mathbb{R}_+$ with respect to the logarithmic metric $d(t,t') = |\ln t- \ln t' |$.

Since $\mathrm{VSO}(\mathbb{R}_+)$ is a proper subset of $L_\infty(\R_+)$, there exist bounded vertical operators which, although being Toeplitz, still do not belong to the $C^*$-algebra generated by Toeplitz operator with bounded vertical symbols. An example of such operator for the unweighted Bergman space case, $\lambda=0$, is given in \cite[Section 6]{HMV}.
Let us recall some details. By \cite[Proposition 6.1]{HMV}, there exists a unique function $\mathbf{f}:\,{\mathbb{R}_+}\to\mathbb{C}$,
such that $\mathbf{f}\in L_1({\mathbb{R}_+},e^{-\eta u}\,du)$ for all $\eta>0$, and the Laplace transform of $\mathbf{f}$:
\[
f(z):=\int_0^{+\infty}\mathbf{f}(u)e^{-zu}\,du,
\]
is given by
\begin{equation*}
f(z)=\dfrac{1}{z+1}\exp\left(\dfrac{i}{3\pi}\ln^2(z+1)\right),
\end{equation*}
where $\ln$ is the principal branch of the natural logarithm function in $\Pi$ (with imaginary part in $(-\pi,\pi]$).
Take now $a(\mathrm{Im}\,z) = \mathbf{f}(2\mathrm{Im}\,z)$. Then by \cite[Proposition 6.2]{HMV}, the corresponding function $\gamma_a$ is given by
\begin{equation*}
\gamma_a(\eta)=\frac{\eta}{\eta+1}\exp\left(\dfrac{i}{3\pi}\ln^2(\eta+1)\right);
\end{equation*}
it belongs to $L_\infty({\mathbb{R}_+})$, but does not belong to $\mathrm{VSO}({\mathbb{R}_+})$.

On the other hand, by \eqref{eq:varphi_delta}, the Toeplitz operator with distributional symbol $\pmb{1}\otimes \delta(y-1)$ does belong to the $C^*$-algebra generated by Toeplitz operator with bounded vertical symbols.

\smallskip
If we drop the condition of the essential boundedness of the function $\varphi$, we obtain a mapping of the space of equivalence classes of $E$-measurable functions to the set of \emph{unbounded} normal vertical operators in $\mathcal{A}^2_{\lambda}(\Pi)$. These operators mutually commute being restricted to the invariant linear set $R^*_{\parr}\Ds(\R_+)$, which is dense in $L_2(\Pi)$.

Finally, note that the vertical operator $T_a$, like any other function of an operator with purely absolutely continuous spectrum,  can never be compact, unless it is zero.

\section{The hyperbolic case, angular operators}

This case has some differences compared with the previous ones. In the first two cases the differential operator $\frac{1}{i}V_{...}$ and the weight factor in the corresponding $L_2$-space involve different coordinates, $\theta$ and $r$ in the first case, and $x$ and $y$ in  the second one. {Now, in the hyperbolic case, the operator has the form $\frac{1}{i}V_{\hyp} = \frac{1}{i} r\frac{\partial }{\partial r}$ while the weight factor equals $r^{\lambda} \sin^{\lambda} \theta$, containing  both the radial coordinate $r$ and the angle $\varphi$ in the weight.} This causes  the operator $\frac{1}{i}V_{\hyp}$ to be not Hermitian in the $\lambda$-dependent space $L_2(\Pi, d\nu_{\lambda})$ anymore. To make it Hermitian we need to correct it by a $\lambda$-dependent non-real additive  term.
It is not difficult to figure out that the required $\lambda$-dependent operator should have the form
\begin{equation*}
 H_{\lambda} = \frac{1}{i}\left(V_{\hyp} + (1 +\lambda/2)\right) = \frac{1}{i} \left( r\frac{\partial }{\partial r} + (1 +\lambda/2)\right).
\end{equation*}
This implies that $H_{\lambda}$ is Hermitian on  $L_2(\Pi, d\nu_{\lambda})$ and can be easily shown to be self-adjoint on the natural domain
\begin{equation*}
\Dc( H_{\lambda})=\left\{f\in L_2(\Pi, d\nu_{\lambda}):\, r\frac{\partial f}{\partial r}\in L_2(\Pi, d\nu_{\lambda})\right\}.
\end{equation*}
This operator commutes with all dilation operators $\pmb{\delta}_\rho: f(z) \mapsto f(\rho^{-1}z), \, \r>0$, providing thus an example of an unbounded angular operator on the half-plane.

We restrict  the operator $H_{\lambda}$ to the Bergman space $\mathcal{A}^2_{\lambda}(\Pi),$ setting $\Mb_\l=H_{\lambda}|_{\mathcal{A}^2_{\lambda}(\Pi)}.$ This operator acts as
  $\Mb_\l= \frac{1}{i}\left(z\frac{\partial }{\partial z} + (1 +\lambda/2)\right).$ It follows from
 \begin{eqnarray*}
 && \left(R_{\hyp}\Mb_{\lambda}R_{\hyp}^* g\right)(\e) = R_{\hyp}\left(\frac{1}{i}\left(z\frac{\partial }{\partial z}+(1+\lambda/2)\right) \right. \\
 && \left. \times  \ \frac{1}{\sqrt{2}} \int_{\mathbb{R}} g(\eta) z^{i\eta-(1+\lambda/2)} \vartheta_{\lambda}(\eta) d\eta\right)
 = (R_{\hyp}  R_{\hyp}^*)(\eta g(\eta)) = \eta g(\eta)
 \end{eqnarray*}
that the operator $\Mb_\l$ is unitarily equivalent to the self-adjoint operator $w(\eta)\mapsto \eta w(\eta)$ in $L_2(\R)$ with its natural domain. Thus, the spectrum of $\Mb_\l$ fills the whole real axis, is absolutely continuous and has multiplicity one.
Further, the spectral measure $E(S)$ of $\Mb_\l$ is obtained by the above unitary transformation from the spectral measure of the multiplication operator, similarly to \eqref{Sp.Measure K},
\begin{equation*}
(E(S)f)(z)=(R_{\hyp}^*\h_S(\eta) R_{\hyp} f)(z),
\end{equation*}
where, as before, $\h_S(\eta)$ is the characteristic function of the Borel set $S\subset \R.$
From the explicit expression for the unitaries $R_{\hyp}, \ R_{\hyp}^*$, it follows that $E(S)$ is an integral operator,
\begin{equation*}
(E(S)f)(z)=\int_{\Pi}\eF_{S} (z,\z)f(\z)d\n_\l(\z),
\end{equation*}
where the kernel $\eF_{S} (z,\z)$ equals
\begin{gather*}
    \eF_{S}(z,\z)=\frac{1}{2}\int_{S}\theta_\l(\eta)^2z^{i\eta-(1+\l/2)}\bar{\zeta}^{-i\eta-(1+\l/2)}d\eta=\\ \nonumber
\frac{1}{2\pi\G(\l+2)}(z\bar{\z})^{-(1+\l/2)}
\int_{S}\left|\G\left(\frac{\l+2}{2}+i\eta\right)\right|^2(z\bar{\z}^{-1})^{i\eta}d\eta.
\end{gather*}
This kernel contains the factor $(z\bar{\z})^{-(1+\l/2)}$  reflecting the main dependence of the weight parameter $\l$ and a universal factor depending on $z\bar{\z}^{-1}.$ This structure of the kernel reflects the dilation invariance of the spectral measure.

A more explicit expression for the resolution of identity $\Eb(\y):=E((-\infty,\y))$, { similar to} the one in  \eqref{EF.2}, is not available.

Following the general scheme of the Spectral Theory, we associate with any function $f\in \Ac^2_\l(\Pi)$ the corresponding scalar measure  $d\r_f(\y)=d\langle f,\Eb(\y) f\rangle.$

We group together now all the statements analogous to those of the previous two cases in one proposition.

\begin{proposition}
 Given a Borel function $\varphi$ on $\mathbb{R}$, the operator
 \begin{equation*}
  \varphi(\Kb) = \int_{\mathbb{R}} \varphi(\eta)\,d\Eb(\eta)
 \end{equation*}
is well defined and normal on its domain
\begin{equation*}
 \Dc_{\varphi} = \left\{f \in \mathcal{A}^2_{\lambda}(\Pi) : \int_{\mathbb{R}} |\varphi(\eta)|^2\, d\rho_f(\eta)  < \infty\right\}.
 \end{equation*}

The operator $\varphi(\Kb)$ is bounded, and thus defined on the whole $\mathcal{A}^2_{\lambda}(\Pi)$, if and only if the function   $\varphi$ is (Lebesgue) essentially bounded.

Defined on such functions, the mapping $\Jc_\l: \varphi\mapsto\varphi(\Mb_\l)$ is an isometric isomorphism of the unital algebra $L_\infty(\R)$ onto a unital commutative algebra in $\Bc(\Ac_\l^2(\Pi)).$ The image of $L_\infty(\R)$ in $\Bc(\Ac_\l^2(\Pi))$ coincides with the von Neumann algebra of all bounded angular operators in $\Ac_\l^2(\Pi).$

Operators $\varphi(\Mb_\l)$ for all Borel functions $\varphi$ form a commutative algebra, with a natural definition of commutativity of unbounded operators, see Sect.3.
\end{proposition}

Let now $T_a$ be an angular Toeplitz operator in $\Ac_\l^2(\Pi)$ with some bounded symbol $a=a(\arg z).$ This operator commutes with dilations and therefore has to be a function of the generating operator $\Mb_\l$, i.e. $T_a=\varphi_a(\Mb_\l)$ for the function $\varphi_a(\eta)$ given by
\begin{equation}\label{fi hyp}
    \varphi_a(\eta)= 2^{\l}(\l+1)\vartheta_\l^2(\eta)\int_0^{\pi} a(\theta)e^{-2\y\theta}\sin^\l\theta d\theta,
\end{equation}
which coincides with $\gamma_a$ in \eqref{gamma hyp}.

Similarly to the previous two cases, the spectral representation of angular operators extends to unbounded radial symbols and to the distributional ones. A detailed analysis of the boundedness conditions in a general setting, similar to the one  made in other cases, see \cite{RozVas},  \cite{RVhp}, is to be published later on. We restrict ourselves to some simple but still enlightening  particular case.

Let $a(\theta)$ be a distribution with compact support in $(0,\pi)$, $a\in \Es'((0,\pi)).$ With this distribution we associate the angular symbol-distribution   $a^\dag = \pmb{1}\otimes a\in C(\R_+)\otimes \Es'((0,\pi)).$  The distribution $\mathsf{a}$ defines the Toeplitz operator $T_{a^\dag}$ in $\Ac_\l^2(\Pi)$ by means of the sesquilinear form
\begin{eqnarray} \label{form.hyp}
 \langle T_{a^\dag}f,g\rangle_{\l} &:=& \textstyle{\frac{2^{\lambda}(\lambda+1)}{\pi}}\pmb{(}a^\dag ,|\Ima z|^{\l}f(z)\overline{g(z)} \pmb{)} \\
 &=& \frac{2^{\lambda}(\lambda+1)}{\pi}\int_{0}^{\infty}r^{\l+1}\pmb{(} a,\sin(\theta)^{\l}f(re^{i\cdot})\overline{g(re^{i\cdot})}\pmb{)}_{\theta}dr, \nonumber
\end{eqnarray}
where $\pmb{(}.,.\pmb{)}_{\theta}$ denotes the action in $\theta$ variable of a distribution in $\Es'((0,\pi))$ on a smooth function on $(0,\pi).$ Note that the expression \eqref{form.hyp} for the sesquilinear form is a natural generalization of the standard expression for the sesquilinear form of a Toeplitz operator in $\Ac_\l^2(\Pi)$ with a bounded symbol.

 The easiest way to study the boundedness of this operator (unlike the much harder case of a general distributional symbol,) is to check that the expression \eqref{gamma hyp} for the spectral function corresponding to this operator produced
  sa bounded function on $\R.$ The natural generalization of the integral in \eqref{gamma hyp} is
 \begin{equation*}
 \pmb{(} a, e^{-2\xi\theta} \sin^{\lambda}\theta \pmb{)}_{\theta}
 \end{equation*}

  Here, again, we use the basic  Theorem  in Sect.4.4 Ch.II.,\cite{G2}, already cited above, to obtain the representation $a=\sum_l D_l b_l$ with functions $b_l\in C_0(0,\pi).$ Similar to the parabolic case, this leads to
\begin{equation}\label{hyp.eta}
  \pmb{(} a, e^{-2\eta\theta} \sin^{\lambda}\theta \pmb{)}_{\theta}=\sum_l\pmb{(} b_l, D_l^*[e^{-2\eta\theta} \sin^{\lambda}\theta ] \pmb{)}_{\theta}.
\end{equation}
{Since  the functions $b_l$ have compact support in $(0,2\pi)$, by \eqref{Gamma}, the corresponding function $\g_{a}(\eta)$ in \eqref{gamma hyp}, calculated using  \eqref{hyp.eta}, is bounded} and, by the spectral representation, this leads to the boundedness of the Toeplitz operator $T_{\mathsf{a}}.$

The explicit description of those $\varphi$ for which
the operator $\varphi(\Kb)$ belongs to the $C^*$-algebra $\mathcal{T}_{\hyp}$ {generated by Toeplitz operator with bounded vertical symbols is given in \cite{EMV}. This happens if and only if $\varphi$ belongs to  $\mathrm{VSO}(\mathbb{R})$, the set of all bounded functions that are uniformly continuous
with respect to the \emph{arcsinh}-metric
\begin{equation*}
 d(t, t') = |\arcsinh t - \arcsinh t'|, \qquad t, t' \in \mathbb{R}.
\end{equation*}
Let us give two examples of unbounded and distributional angular symbols, for which the  corresponding Toeplitz operators belong nevertheless to the $C^*$-algebra $\mathcal{T}_{\hyp}$.

First, let $a(\theta) = \theta^{-\beta} \sin \theta^{-\alpha}$, with $0 < \beta < 1$ and $\alpha > 0$. This symbol, oscillating near $0$ and unbounded,  was considered in \cite[Example 8.5.3]{book}; as shown there, the corresponding function $\gamma_a(\eta)$ is continuous on $\R$ and has finite limits as $\eta \to \pm \infty$. Therefore, the Toeplitz operator $T_a$ belongs to $\mathcal{T}_{\hyp}$.

For the second example, consider the distributional symbol $a^\dag=\pmb{1}\otimes \delta(\theta - \frac{\pi}{2})$. Then, by \eqref{fi hyp} and \eqref{eq:vartheta},
\begin{equation*}
 \varphi_{a^\dag}(\eta) = 2^{\l}(\l+1)\vartheta_\l^2(\eta)\pmb{(} \delta(\theta - \textstyle{\frac{\pi}{2}}), e^{-2\eta\theta} \sin^{\lambda}\theta \pmb{)}_{\theta} =
 \displaystyle{\frac{2^{\lambda}}{\pi \G(\lambda+1)}}\left|\G(\textstyle{\frac{\lambda+1}{2}} + i \eta)\right|^2.
\end{equation*}
Function $\varphi_{a^\dag}$ is continuous on $\R$ and, by \cite[Formula 8.328.1]{GR}, $\varphi_{a^\dag}(\eta) \to 0$ as $\eta \to \pm \infty$, which implies that the Toeplitz operator $T_{a^\dag}$ is bounded and belongs to
$\mathcal{T}_{\hyp}$.

Similarly to the parabolic case, a nonzero angular operator, being a function of an operator with absolutely continuous spectrum, can never be compact.

\end{document}